\documentclass[12pt,a4paper,twoside]{article}
\usepackage{mathtext}
\usepackage[cp1251]{inputenc}
\usepackage[T2A]{fontenc}
\usepackage[russian]{babel}
\usepackage[dvips]{graphicx}
\usepackage{amsmath}
\usepackage{amssymb}
\usepackage{amsxtra}
\usepackage{latexsym}
\usepackage{ifthen}
\usepackage[centerlast,footnotesize,sl]{caption2}
\begin{document}

\newcounter{lemma}[section]
\newcommand{\lemma}{\par \refstepcounter{lemma}%
{\bf Лемма \arabic{section}.\arabic{lemma}.}}
\renewcommand{\thelemma}{\thesection.\arabic{lemma}}

\newcounter{corollary}[section]
\newcommand{\corollary}{\par \refstepcounter{corollary}%
{\bf Следствие \arabic{section}.\arabic{corollary}.}}
\renewcommand{\thecorollary}{\thesection.\arabic{corollary}}

\newcounter{remark}[section]
\newcommand{\remark}{\par \refstepcounter{remark}%
{\bf Замечание \arabic{section}.\arabic{remark}.}}
\renewcommand{\theremark}{\thesection.\arabic{remark}}

\newcounter{theorem}[section]
\newcommand{\theorem}{\par \refstepcounter{theorem}%
{\bf Теорема \arabic{section}.\arabic{theorem}.}}
\renewcommand{\thetheorem}{\thesection.\arabic{theorem}}

\newcounter{proposition}[section]
\newcommand{\proposition}{\par \refstepcounter{proposition}%
{\bf Предложение \arabic{section}.\arabic{proposition}.}}
\renewcommand{\theproposition}{\thesection.\arabic{proposition}}

\newcommand{\proof}{{\it Доказательство.\,\,}}

\renewcommand{\theequation}{\arabic{section}.\arabic{equation}}

\def\Xint#1{\mathchoice
   {\XXint\displaystyle\textstyle{#1}}
   {\XXint\textstyle\scriptstyle{#1}}
   {\XXint\scriptstyle\scriptscriptstyle{#1}}
   {\XXint\scriptscriptstyle\scriptscriptstyle{#1}}
   \!\int}
\def\XXint#1#2#3{{\setbox0=\hbox{$#1{#2#3}{\int}$}
     \vcenter{\hbox{$#2#3$}}\kern-.5\wd0}}
\def\dashint{\Xint-}

\noindent
\renewcommand{\baselinestretch}{1.25}
\normalsize \large

\medskip

\centerline{\bf On Neumann and Poincare problems in mathematical
physics}

\medskip

\centerline{\bf Artyem Yefimushkin}

\medskip

\centerline{ Inst. Appl. Math. Mech., National Academy of Sciences
of Ukraine}

\centerline{E-mail: a.yefimushkin@gmail.com}

\medskip

\centerline{\bf Abstract}

It is proved the existence of nonclassical solutions of the Neumann
and Poincare problems for generalizations of the Laplace equation in
anisotropic and nonhomogeneous media in almost smooth domains with
arbitrary boundary data that are measureable with respect to
logarithmic capacity. Moreover, it is shown that the spaces of such
solutions have the infinite dimension.

\medskip


\bigskip

{\bf 1. Introduction} The classical boundary value problems of the
theory of analytic functions, such as the Dirichlet, Hilbert,
Riemann, Neumann and Poincare problems, take a fundamental part in
contemporary analysis and applications to actual problems of
mathematical physics. Note that these boun\-da\-ry value problems
are closely interconnected, see e.g. \cite{AIM}--\cite{UMV}. In this
paper we continue the development of the theory of the above
boundary value problems for analytic functions in Jordan domains
with mea\-su\-rab\-le boundary data and extend the theory to the
more general class of quasiconformal functions.

It is well-known that the Neumann problem for Laplace equation has
no classical solutions generally speaking even for some continuous
boundary data, see e.g. \cite{Mi}. The main goal of this paper is to
show that the Neumann problem with arbitrary boundary data that are
measurable with respect to logarithmic capacity has nonclassical
solutions for the Laplace equation as well as for its
ge\-ne\-ra\-li\-za\-tions. The result is based on a reduction of
this problem to the Riemann-Hilbert boundary value problem whose
solution was recently obtained in \cite{UMV}.

First of all recall a more general {\bf problem on directional
derivatives} for the harmonic functions in the unit disk $\mathbb D
= \{ z\in\mathbb{C}: |z|<1\}$, $z=x+iy$. The classic setting of the
latter problem is to find a function $u:\mathbb D\to\mathbb R$ that
is twice continuously differentiable, admits a continuous extension
to the boundary of $\mathbb D$ together with its first partial
derivatives, satisfies the Laplace equation
$$
\Delta u\ :=\ \frac{\partial^2 u}{\partial x^2}\ +\ \frac{\partial^2
u}{\partial y^2}\ =\ 0 \quad\quad\quad\forall\ z\in\mathbb D $$ and
the boundary condition with a prescribed continuous date $\varphi :
\partial\mathbb D\to\mathbb R$:
$$
\frac{\partial u}{\partial \nu}\ =\ \varphi(\zeta) \quad\quad\quad
\forall\ \zeta\in\partial\mathbb D
$$
where $\frac{\partial u}{\partial \nu}$ denotes the derivative of
$u$ at $\zeta$ in a direction $\nu = \nu(\zeta)$, $|\nu(\zeta)|=1$:
$$
\frac{\partial u}{\partial \nu}\ :=\ \lim_{t\to 0}\
\frac{u(\zeta+t\cdot\nu)-u(\zeta)}{t}\ .
$$

{\bf The Neumann problem} is a special case of the above problem on
directional derivatives with  the boundary condition
$$
\frac{\partial u}{\partial n}\ =\ \varphi(\zeta) \quad\quad\quad
\forall\ \zeta\in\partial\mathbb D
$$
where $n$ denotes the unit interior normal to $\partial\mathbb D$ at
the point $\zeta$.


In turn, the above problem on directional derivatives is a special
case of {\bf the Poincare problem} with  the boundary condition
$$
a\cdot u\ +\ b\cdot\frac{\partial u}{\partial \nu}\ =\
\varphi(\zeta) \quad\quad\quad \forall\ \zeta\in\partial\mathbb D
$$
where $a=a(\zeta)$ and $b=b(\zeta)$ are real-valued functions given
on $\partial\mathbb D$.


Recall also that twice continuously differentiable solutions of the
Laplace equation are called {\bf harmonic functions}. As known, such
functions are infinitely differentiable as well as they are real and
imaginary parts of analytic functions.

\bigskip

{\bf 2. Definitions and preliminary remarks.} The partial
differential equa\-tions in the divergence form below take a
sig\-ni\-fi\-cant part in many problems of mathematical physics, in
particular, in anisotropic and inhomogeneous media. These equations
are closely interconnected with Beltrami equations, see e.g.
\cite{AIM}.

Let $D$ be a domain in the complex plane $\mathbb{C}$ and let
$\mu:D\to\mathbb{C}$ be a measurable function with $|\mu(z)|<1$ a.e.
Recall that a partial differential equation
$$ f_{\bar{z}}=\mu(z)\cdot f_z\
 \eqno(1)$$
where $f_{\bar z}={\bar\partial}f=(f_x+if_y)/2$, $f_{z}=\partial
f=(f_x-if_y)/2$, $z=x+iy$, $f_x$ and $f_y$ are partial derivatives
of the function $f$ in $x$ and $y$, respectively, is said to be a
{\bf Beltrami equation}. The Beltrami equation (1) is said to be
{\bf nondegenerate} if $||\mu||_{\infty}<1$.

In this connection, note that if $f=u+i\cdot v$ is a regular
solution of the Beltrami equation (1), then the function $u$ is a
continuous generalized solution of the divergence-type equation
$$
{\rm div}\, A(z)\nabla\,u=0\ ,
 \eqno(2)$$
called {\it A-harmonic function}, i.e. $u\in C\cap W^{1,1}$ and
$$
\int\limits_D \langle A(z)\nabla
u,\nabla\varphi\rangle=0\,\,\,\,\,\,\,\,\,\,\,\,\forall\ \varphi\in
C_0^\infty(D)\ ,
$$
where $A(z)$ is the matrix function: $$
A=\left(\begin{array}{ccc} {|1-\mu|^2\over 1-|\mu|^2}  & {-2{\rm Im}\,\mu\over 1-|\mu|^2} \\
                            {-2{\rm Im}\,\mu\over 1-|\mu|^2}          & {|1+\mu|^2\over 1-|\mu|^2}  \end{array}\right).
\eqno(3)$$.

As we see in (3), the matrix $A(z)$ is symmetric and its entries
$a_{ij}=a_{ij}(z)$ are dominated by the quantity
$$
K_{\mu}(z)\ =\ \frac{1\ +\ |\mu(z)|}{1\ -\ |\mu(z)|}\ ,
$$
and, thus, they are bounded if the Beltrami equation (1) is not
degenerate.

Vice verse, uniformly elliptic equations (2) with symmetric $A(z)$
and ${\rm det}\,A(z)\equiv 1$ just correspond to nondegenerate
Beltrami equations (1) with coefficient
$$
\mu\ =\ \frac{1}{\mathrm{det}\, (I+A)}\ (a_{22}-a_{11}\ -\
2ia_{21})\ =\ \frac{a_{22}-a_{11}\ -\ 2ia_{21}}{1\ +\ \mathrm{Tr}\,
A\ +\ \mathrm{det}\, A}\ ,\eqno(4)
$$
where $I$ denotes identity $2\times 2$ matrix $\mathrm{Tr}\, A =
a_{22}+a_{11}$, see e.g. theorem 16.1.6 in \cite{AIM}. Following
\cite{GRY}, call all such matrix functions $A(z)$ of the {\it class
${\cal{B}}$}. Recall that equation (2) is said to be {\bf uniformly
elliptic}, if $a_{ij}\in L^{\infty}$ and $\langle A(z)\eta,\eta
\rangle\geq\varepsilon|\eta|^2$ for some $\varepsilon>0$ and for all
$\eta\in\mathbb{R}^2$.

Finally, recall that homeomorphic solutions of Beltrami equations
(1) of class $W^{1,1}_{\rm loc}$ is said to be {\bf quasiconformal
mappings}, see e.g. \cite{Alf,LV}. The images of the unit disk
$\mathbb D = \{ z\in\mathbb{C}: |z|<1\}$ under quasiconformal
mappings of $\mathbb{C}$ onto itself are called {\bf quasidisks},
and their boundaries are called {\bf quasicircles} or {\bf
quasiconformal curves}. Recall that the bijective continuous image
of a circle in $\mathbb{C}$ is called a {\bf Jordan curve}. As
known, any smooth (or Lipschitzian) Jordan curve is a quasiconformal
curve, see e.g. point II.8.10 in \cite{LV}.

\bigskip

{\bf 3. On Neumann and Poincare problems for harmonic functions}

\bigskip

Let us start with the unit disk because proofs in this case are more
direct and clear.

{\bf Theorem 1.} {\it Let $\nu:\partial\mathbb D\to\mathbb{C}$,
$|\nu(\zeta)|\equiv1$ be a function of bounded variation, and let
$\varphi:\partial\mathbb D\to\mathbb{R}$ be a measurable function
with respect to logarithmic capacity. Then there exist harmonic
functions $u:\mathbb D\to\mathbb{R}$ such that $$ \lim_{z\to\zeta}\
\frac{\partial u}{\partial\nu}\ =\ \varphi(\zeta) \eqno(5)$$ along
any nontangential paths for a.e. $\zeta\in\partial\mathbb D$ with
respect to logarithmic capacity.}

{\it Proof.} Indeed, by Proposition 6.1 in \cite{UMV}, there exist
analytic function $f:D\to\mathbb{C}$ such that
$$\lim_{z\to\zeta}\ \mathrm{Re}\ \overline{\nu(\zeta)\cdot f(z)}\ =\
\lim_{z\to\zeta}\ \mathrm{Re}\ {\nu(\zeta)}\cdot f(z)\ =\
{\varphi(z)} \eqno(6)$$ along any nontangential paths for a.e.
$\zeta\in\partial \mathbb D$ with respect to logarithmic capacity.
Note that an indefinite integral $F$ of $f$ in $\mathbb{D}$ is also
an analytic function and, correspondingly, the harmonic functions
$u=\mathrm{Re}\, F$ and $v=\mathrm{Im}\, F$ satisfy the system of
Cauchy-Riemann $v_x=-u_y$ и $v_y=u_x$. Hence $$f\ =\ F'\ =\ F_x\ =\
u_x\ +\ i\cdot v_x\ =\ u_x\ -\ i\cdot u_y\ =\ \overline{\nabla u}$$
where $\nabla u=u_x+i\cdot u_y$ is the gradient of the function $u$
in complex form. Thus, (5) follows from (6), i.e. $u$ is one of the
desired harmonic functions because its directional derivative
$$\frac{\partial u}{\partial\nu}\ =\
{\mathrm Re}\, \overline{\nu}\cdot\nabla u\ =\ {\mathrm Re}\,
{\nu}\cdot\overline{\nabla u}\ =\ \langle \nu,\nabla u\rangle$$ is
the scalar product of $\nu$ and the gradient $\nabla u$ interpreted
as vectors in $\mathbb{R}^2$. $\Box$

\medskip

{\bf Remark 1.} We are able to say more in the case $\mathrm {Re}\
n\cdot\overline{\nu}>0$ where $n=n(\zeta)$ is the unit interior
normal with a tangent to $\partial\mathbb{D}$ at the point
$\zeta\in\partial\mathbb{D}$. In view of (5), since the limit
$\varphi(\zeta)$ is finite, there is a finite limit $u(\zeta)$ of
$u(z)$ as $z\to\zeta$ in $\mathbb{D}$ along the straight line
passing through the point $\zeta$ and being parallel to the vector
$\nu(\zeta)$ because along this line, for $z$ and $z_0$ that are
close enough to $\zeta$, $$ u(z)\ =\ u(z_0)\ -\
\int\limits_{0}\limits^{1}\ \frac{\partial u}{\partial \nu}\
(z_0+\tau (z-z_0))\ d\tau\ .$$ Thus, at each point with the
condition (5), there is the directional derivative
$$
\frac{\partial u}{\partial \nu}\ (\zeta)\ :=\ \lim_{t\to 0}\
\frac{u(\zeta+t\cdot\nu)-u(\zeta)}{t}\ =\ \varphi(\zeta)\ .
$$

\medskip

In particular, in the case of the Neumann problem, we have by
Theorem 1 and Remark 1 the following significant result.

\medskip

{\bf Theorem 2.} {\it For each function $\varphi:\partial\mathbb
D\to\mathbb R$ that is measurable with respect to logarithmic
capacity, one can find harmonic functions $u:\mathbb D\to\mathbb C$
such that, for a.e. point $\zeta\in\partial\mathbb D$ with respect
to logarithmic capacity, there exist:

\bigskip

1) the finite radial limit
$$
u(\zeta)\ :=\ \lim\limits_{r\to 1}\ u(r\zeta)$$

2) the normal derivative
$$
\frac{\partial u}{\partial n}\, (\zeta)\ :=\ \lim_{t\to 0}\
\frac{u(\zeta+t\cdot n)-u(\zeta)}{t}\ =\ \varphi(\zeta)
$$

3) the nontangential limit $$ \lim_{z\to\zeta}\ \frac{\partial
u}{\partial n}\, (z)\ =\ \frac{\partial f}{\partial n}\, (\zeta)$$
where $n=n(\zeta)$ denotes the unit interior normal to
$\partial\mathbb D$ at the point $\zeta$.}

\bigskip

Recall that a Jordan domain is called {\bf Lipschitzian} if its
boundary is bilipschitzian image of a circle. It is clear that such
curve is rectifiable, and rectifiable curves have tangent to almost
all points with respect to the length measure. A Jordan domain is
called {\bf almost smooth} if it is Lipschitzian and has tangent to
almost all points with respect to logarithmic capacity.

\medskip

{\bf Theorem 3.} {\it Let $D$ be an almost smooth Jordan domain in
the complex plane $\mathbb{C}$, $\nu:\partial D\to~\mathbb{C}$,
$|\nu(\zeta)|\equiv1$ be a function of bounded variation and let
$\varphi:\partial D\to\mathbb{R}$ be a function that is measurable
with respect to logarithmic capacity. Then there exist harmonic
functions $u:D\to\mathbb{R}$ such that $$\lim_{z\to\zeta}\
\frac{\partial u}{\partial\nu}\ =\ \varphi(\zeta) \eqno(7)$$ along
any nontangential paths for a.e. point $\zeta\in\partial D$ with
respect to logarithmic capacity. }

\medskip

{\it Proof.} The case of almost smooth Jordan domains $D$ reduces to
the case of the unit disk $\mathbb{D}$ in the following way. First
of all, by the Riemann theorem one can find conformal mapping
$\omega$ of the domain $D$ onto $\mathbb{D}$, see e.g. Theorem
II.2.1 in \cite{Go}. Then by Caratheodory theorem $\omega$ is
extended to a homeomorphism of $\overline{D}$ onto
$\overline{\mathbb{D}}$, see e.g. Theorem II.C.1 in \cite{Ku}.

As it was noted in Section 2, the boundaries of Lipschitzian domains
are quasiconformal curves. Thus, by the reflection principle for
quasiconformal mappings, involving the conformal reflection
(inversion) with respect to the unit circle in the image and a
quasiconformal reflection with respect to $\partial D$ in the
preimage, we are able to extend $\omega$ to a quasiconformal mapping
$\Omega:\mathbb{C}\to\mathbb{C}$, see e.g. I.8.4, II.8.2 and II.8.3
in \cite{LV}. It is clear also that ${\cal{N}} :=\nu\circ \Omega
^{-1}|_{\partial\mathbb D}$ is a function of bounded variation,
$V_{\cal{N}}(\partial{\mathbb{D}})=V_{\nu}(\partial{{D}})$.

The logarithmic capacity of a set coincides with its transfinite
diameter, see e.g.  \cite{F} and the point 110 in \cite{N}.
Moreover, quasiconformal mappings are H\"{o}lder continuous on
compacta, see e.g. Theorem II.4.3 in \cite{LV}. Hence the mappings
$\Omega$ and $\Omega^{-1}$ transform sets of logarithmic capacity
zero on $\partial D$ into sets of logarithmic capacity zero on
$\partial \mathbb{D}$ and vice versa.

The function $\Phi :=\varphi\circ \Omega ^{-1}|_{\partial\mathbb D}$
is measurable with respect to logarithmic capacity. Indeed, under
the given mappings any sets that are measurable with respect to
logarithmic capacity are transformed into sets that are measurable
with respect to logarithmic capacity because any such a set can be
represented in the form of the union of a sigma-compact set and a
set of logarithmic capacity zero and, under continuous mappings, the
compact sets are transformed into compact sets and the latter sets
are measurable with respect to logarithmic capacity.

Next, by the Lindel\"{o}f lemma $\arg\
[\omega(\zeta)-\omega(z)]-\arg\ [\zeta-z]\to\mathrm const$ as
$z\to\zeta$ for every point $\zeta\in\partial D$ with a tangent to
$\partial D$, see e.g. Theorem II.C.2 in \cite{Ku}. Hence
nontangential paths in $D$ to $\zeta\in\partial D$ are transformed
under the mapping $\omega$ into nontangential paths in $\mathbb D$
to $\xi = \omega(\zeta)\in\partial\mathbb D$ for a.e.
$\zeta\in\partial D$ with respect to logarithmic capacity. And vice
versa, nontangential paths in $\mathbb D$ to $\xi\in\partial\mathbb
D$ are transformed under the mapping $\omega^{-1}$ into
nontangential paths in $D$ to $\zeta = \omega^{-1}(\xi)\in\partial
D$ for a.e. $\xi\in\partial\mathbb D$ with respect to logarithmic
capacity.

By the Theorem 1 one can find a harmonic function $U:\mathbb
D\to\mathbb{R}$ such that $$\lim_{w\to\xi}\ \frac{\partial
U}{\partial{\cal{N}}}\ (w)\ =\ \Phi(\xi) \eqno(8)$$ along any
nontangential paths to a.e. point $\xi\in\partial\mathbb D$ with
respect to logarithmic capacity. Moreover, in the simply connected
domain $\mathbb{D}$, one can find a harmonic function $V:\mathbb
D\to\mathbb R$ such that $g:=U+i\cdot V$ is a single-valued analytic
function in $\mathbb{D}$.

Let $F$ be an indefinite integral of the analytic function
$g^{\prime}\cdot(\omega^{-1})^{\prime}$ in $\mathbb D$ and let
$f:=F\circ \omega$. Then $F$ and $f$ are also single-valued analytic
functions in $\mathbb{D}$ and $D$, correspondingly, and elementary
calculations show that $$ f^{\prime}\ =\
F^{\prime}\circ\omega\cdot\omega^{\prime}\ =\
F^{\prime}\circ\omega\cdot(\omega^{\prime}\circ\omega^{-1})\circ\omega\
=\ [F^{\prime}/(\omega^{-1})^{\prime}]\circ\omega\ =\
g^{\prime}\circ\omega\ .
$$
Thus,
$$
\frac{\partial f}{\partial{\nu}}\ =\ \nu\cdot f^{\prime}\ =\
\nu\cdot g^{\prime}\circ\omega\ =\ \frac{\partial
g}{\partial{\cal{N}}}\circ\omega\ ,
$$
where $\nu = \nu(\zeta)$, $\zeta\in\partial D$, and
${\cal{N}}={\cal{N}}(\xi)$, $\xi = \omega(\zeta)\in\partial\mathbb
D$. Hence, for $u:={\mathrm Re}\, f$, we have the equality
$$
\frac{\partial u}{\partial{\nu}}\ =\ \frac{\partial
U}{\partial{\cal{N}}}\circ\omega
$$
and, consequently, $u$ is a desired harmonic function by (8). $\Box$

\bigskip

By Theorem 3, arguing similarly to Remark 1, we have the following
sig\-ni\-fi\-cant result on the Neumann problem.

\bigskip

{\bf Theorem 4.} {\it Let $D$ be an almost smooth Jordan domain in
the complex plane $\mathbb{C}$ and let a function $\varphi:\partial
D\to\mathbb{R}$ be measurable with respect to logarithmic capacity.
Then one can find harmonic functions $u: D\to\mathbb{C}$ such that
for a.e. $\zeta\in\partial D$ with respect to logarithmic capacity,
there exist:

\bigskip

1) the finite normal limit
$$
u(\zeta)\ :=\ \lim_{z\to\zeta}\ u(z)$$

2) the normal derivative
$$
\frac{\partial u}{\partial n}\, (\zeta)\ :=\ \lim_{t\to0}\
\frac{u(\zeta+t\cdot n)\ -\ u(\zeta)}{t}\ =\ \varphi(\zeta)
$$

3) the nontangential limit $$ \lim_{z\to\zeta}\ \frac{\partial
u}{\partial n}\, (z)\ =\ \frac{\partial u}{\partial n}\, (\zeta)\
.$$}

\bigskip

{\bf 4. On Neumann and Poincare problems for $A$-harmonic functions}

\bigskip

{\bf Theorem 5.} {\it Let $D$ be an almost smooth Jordan domain in
$\mathbb{C}$, $A(z), \:z\in D$, be a matrix function of the class
$\mathcal{B}\cap C^{\alpha},\:\alpha\in(0,1)$, $\nu:\partial
D\to~\mathbb{C}$, $|\nu(\zeta)|\equiv1$, be a function of bounded
variation and let a function $\varphi:\partial D\to\mathbb{R}$ be
measurable with respect to logarithmic capacity. Then there exist
$A$-harmonic functions $u:D\to\mathbb{R}$ of the class
$C^{1+\alpha}$ such that $$\lim_{z\to\zeta}\ \frac{\partial
u}{\partial\nu}\, (z)\ =\ \varphi(\zeta) \eqno(9)$$ along any
nontangential paths for a.e. $\zeta\in\partial D$ with respect to
logarithmic capacity.}

{\it Proof.} By remarks in Section 2, a desired function $u$ is a
real part of a solution $f$ of the class $W^{1,1}_{\mathrm loc}$ for
corresponding Beltrami equation with $\mu\in C^{\alpha}$, see e.g.
Theorem 16.1.6 in \cite{AIM}. By Lemma 1 in \cite{GRY} $\mu$ is
extended to a H\"{o}lder continuous function
$\mu_*:\mathbb{C}\to\mathbb{C}$ of the class  $C^{\alpha}$. Hence
also, for every $k_*\in(k,1)$, there is an open neighborhood $U$ of
$D_*$ such that $|\mu(z)|<k_*$. Let $D_*$ be a connected component
of $U$ containing $\overline D$.

Next, there is a quasiconformal mapping $h:{{D_*}}\to{\mathbb{C}}$
a.e. satisfying the Beltrami equation (1) with the complex
coefficient $\mu^*:=\mu_*|_{D_*}$ in $D_*$, see e.g. Theorem V.B.3
in \cite{Alf}. Note that the mapping $h$ has the H\"older continuous
first partial derivatives in $D_*$ with the same order of the
H\"older continuity as $\mu$, see e.g. \cite{Iw} and also
\cite{IwDis}. Moreover, the mapping $h$ is regular, i.e. its
Jacobian $$ J_h(z)\ne 0\  \ \ \ \ \ \ \ \ \ \ \ \ \forall\ z\in D_*\
, \eqno(10)$$ see e.g. Theorem V.7.1 in \cite{LV}. Thus, the
directional derivative $$h_{\omega}(z)=\frac{\partial
h}{\partial\omega}\, (z)\ :=\ \lim_{t\to0}\ \frac{h(z\ +\
t\cdot\omega)\ -\ h(z)}{t}\ \ne\ 0 \ \ \ \ \ \ \ \ \ \ \ \ \forall\
z\in D_*\ \forall\ \omega\in\partial\mathbb D$$ and it is continuous
by the collection of the variables $\omega\in\partial\mathbb D$ and
$z\in D_*$. Thus, the functions
$$ \nu_*(\zeta)\ :=\
\frac{|h_{\nu(\zeta)}(\zeta)|}{h_{\nu(\zeta)}(\zeta)}\ \ \ \ and \ \
\ \ \varphi_*(\zeta)\ :=\ \frac{\varphi(\zeta)}
{|h_{\nu(\zeta)}(\zeta)|}
$$ are measurable with respect to logarithmic capacity, see e.g.
17.1 in \cite{KZPS}.

The logarithmic capacity of a set coincides with its transfinite
diameter, see e.g.  \cite{F} and the point 110 in \cite{N}.
Moreover, quasiconformal mappings are H\"{o}lder continuous on
compacta, see e.g. Theorem II.4.3 in \cite{LV}. Hence the mappings
$h$ and $h^{-1}$ transform sets of logarithmic capacity zero on
$\partial D$ into sets of logarithmic capacity zero on $\partial
D^*$, where $D^*:=h(D)$, and vice versa.

Further, the functions ${\cal{N}} := \nu_*\circ h^{-1}|_{\partial
D^*}$ and $\Phi :=\left(\varphi_* /h_{\nu}\right)\circ
h^{-1}|_{\partial D^*}$ are mea\-su\-rab\-le with respect to
logarithmic capacity. Indeed, a measurable set with respect to
lo\-ga\-rith\-mic capacity is transformed under the mappings $h$ and
$h^{-1}$ into measurable sets with respect to logarithmic capacity
because such a set can be represented as the union of a
sigma-com\-pac\-tum and a set of logarithmic capacity zero and
compacta are transformed  under continuous mappings into compacta
and compacta are measurable with respect to logarithmic capacity.

Recall that the distortion of angles under quasiconformal mappings
$h$ и $h^{-1}$ is bounded, see e.g. \cite{A}-\cite{Ta}. Thus,
nontangential paths to $\partial D$ are transformed into
nontangential paths to $\partial D^*$ for a.e. $\zeta\in\partial D$
with respect to logarithmic capacity and inversly.

By Theorem 3, one can find a harmonic function $U:D^*\to\mathbb{R}$
such that $$\lim_{w\to\xi}\ \frac{\partial U}{\partial\mathcal{N}}\,
(w)\ =\ \Phi(\xi) \eqno(11)$$ along any nontangential paths for a.e.
$\xi\in\partial D^*$ with respect to logarithmic capacity.

Moreover, one can find a harmonic function $V$ in the simply
connected domain $D^*$ such that $F=U+iV$ is an analytic function
and, thus, $u:=\mathrm{Re}\, f=U\circ h$, where $f:=F\circ h$, is a
desired $A$-harmonic function because $f$ is a regular solution of
the corresponding Beltrami equation (1) and also $$ u_{\nu} =
\langle\ \nabla U \circ h\ ,\ h_{\nu}\ \rangle = \langle\
\nu_*\cdot\nabla U \circ h\ ,\ \nu_*\cdot \, h_{\nu}\ \rangle =$$
$$=\langle\ \frac{\partial U}{\partial\mathcal{N}}\ \circ h\ ,\
\nu_*\cdot\, h_{\nu}\ \rangle = \frac{\partial
U}{\partial\mathcal{N}} \ \circ h\ \cdot\ {\mathrm Re}\,
(\nu_*h_{\nu}),
$$ i.e. condition (9) holds along any nontangential paths for a.e. $\zeta\in\partial D$ with respect
to logarithmic capacity.  $\Box$

\bigskip

The following statement concerning to the Neumann problem for
$A$-harmonic functions is a special significant case of Theorem 5.

\medskip

{\bf Theorem 6.} {\it Let $D$ be an almost smooth Jordan domain in
$\mathbb{C}$, the interior unit normal $n=n(\zeta)$ to $\partial{D}$
has bounded variation, $A(z), \:z\in D$, be a matrix function of
class $\mathcal{B}\cap C^{\alpha},\:\alpha\in(0,1)$ and let a
function $\varphi:\partial D\to\mathbb{R}$ be measurable with
respect to logarithmic capacity. Then there exist $A$-harmonic
function $u:D\to\mathbb{R}$ of class $C^{1+\alpha}$ such that for
a.e. $\zeta\in\partial D$ with respect to logarithmic capacity there
exist:

\bigskip

1) the finite normal limit
$$
u(\zeta)\ :=\ \lim_{z\to\zeta}\ u(z)$$

2) the normal derivative
$$
\frac{\partial u}{\partial n}\, (\zeta)\ :=\ \lim_{t\to0}\
\frac{u(\zeta+t\cdot n)\ -\ u(\zeta)}{t}\ =\ \varphi(\zeta)
$$

3) the nontangential limit $$ \lim_{z\to\zeta}\ \frac{\partial
u}{\partial n}\, (z)\ =\ \frac{\partial u}{\partial n}\, (\zeta)\
.$$}

In particular, in the unit disk $\mathbb{D}$, the unit interior
normal $n=n(\zeta)$ to $\partial\mathbb{D}$ has bounded variation
and, thus, the conclusions 1-3 of the latter theorem hold.

\medskip

{\bf 5. On the dimension of spaces of solutions}

\medskip

{\bf Theorem 7.} {\it The spaces of solutions in Theorems 1-6 have
the infinite dimension.}

{\it Proof.} In view of the equivalence of the problem on the
directional derivatives to the corresponding Riemann-Hilbert
boundary value problem established under the proof of Theorem 1, the
conclusion of Theorem 7 follows from Theorem 8.2 in \cite{UMV}
because Theorems 2-6 are successively reduced to Theorem 1. $\Box$

\medskip

{\bf Acknowledgement.} I would like to thank Professor Vladimir
Ryazanov for his scientific supervision.

\medskip

\end{document}